\documentclass[12pt]{article}
\usepackage{graphics,psfrag,amsbsy}
\usepackage[mathscr]{eucal}
\usepackage{euler}
\usepackage{amsmath}
\usepackage{amssymb}
\usepackage{natbib}

\def\email#1{\date{\tt#1}}
\def\address#1{\par\noindent#1\smallskip}

\begin{document}

\title{Steady-State Creep Analysis of Pressurized Pipe Weldments by Perturbation Method}

\author{A. Shutov, H. Altenbach, K. Naumenko}

\email{shutov@ngs.ru, Holm.Altenbach@iw.uni-halle.de, Konstantin.Naumenko@iw.uni-halle.de}

\maketitle

\begin{abstract}

The stress analysis of pressurized circumferential pipe weldments
under steady state creep is considered. The creep response of the
material is governed by Norton's law. Numerical and analytical
solutions are obtained by means of perturbation method, the
unperturbed solution corresponds to the stress field in a
homogeneous pipe. The correction terms are treated as stresses
defined with the help of an auxiliary linear elastic problem.
Exact expressions for jumps
of hoop and radial stresses at the interface are obtained.
The proposed technique essentially simplifies parametric analysis of
multi-material components.

\end{abstract}



\section{Introduction}

The welded pipelines subjected to high pressure and temperature
are widely used in different branches of industry. Under such
conditions, the creep and damage effects should be taken into
account for accurate assurance of long-term reliability
\citep{Roche92}. Application of computational continuum creep damage
mechanics (see, for example \citet{Altenbach2001}, \citet{Hayhurst2001a})
coupled with increasing power of computers can
accomplish this task.
In recent years the finite element method has become the widely accepted
tool for the structural analysis in the creep range \citep{Hayhurst2001a}. A
user defined creep material subroutine with appropriate constitutive and
evolution equations can be developed and incorporated into the commercial
finite element code to perform a numerical time step solution of creep and
long term strength problems.

On the other hand, in addition to more and
more sophisticated numerical analysis, simplified models of creep
response are required. These models should provide a better
intuitive insight into the problem and give a quantitative
description of the solution.

The assessment of reliability of user-defined creep material subroutines and the choice
of suitable numerical parameters like the element type, the mesh density, and time step control
are complicated problems, particulary if studying creep of multi-material structures.
Therefore, it is important to have reference solutions of benchmark
problems. Such solutions should be obtained by use of
alternative analytical or semi-analytical methods which do not require the
spatial discretization techniques and allow for studying the behavior of
stress and deformation gradients.
The objective of this paper is to develop an alternative semi-analytical
solutions to creep problems for multi-material pipe structures. Particularly
we address the analysis of stress gradients in the local zones of
material connections.

To obtain a semi-analytical solution we shall make the following
simplifications. We assume the idealized material behavior having the secondary creep stage
only. In this case the steady state solution of creep in the pipe exists,
for which the stresses do not depend on time. We assume that the difference between the material properties of
constituents is not great. Particularly the difference between the minimum
creep rates for the same stress level should not exceed the value of 2.

The lifetime of a welded pipe under creep conditions is less than
that of homogeneous one. The effect of reliability reduction is of
big interest, therefore large numbers of model problems were
proposed. The most commonly used approach simulates weldment as a
region with non-uniformly distributed material properties
(\citet{Browne1981}, \citet{Coleman1985}, \citet{Hall1991},
\citet{Perrin1999}, \citet{Hayhurst2001b}, \citet{Hyde2003}).

Within the framework of this approach it is often necessary to
consider a number of parameter distribution cases. The amount of
problems to be solved increases with the
number of changing material parameters and parametric analysis
becomes very complicated.

 Drawing an analogy with some simple systems,
for which an analytical solution is available (\citet{Hyde1996,Hyde2000},
\citet{Naumenko2005}) is useful for understanding how the
parameter change can affect the solution. However this is not
enough for a proper estimation of stress distribution.

Two main types of constitutive equations are often used in weld
modelling, namely, Norton's steady-state creep law and continuum
damage mechanics equations for tertiary creep \citep{Kachanov1986}.
Below only the constitutive equations of Norton's law are
considered. Even by this assumption the structural response may be captured
very well. Some special techniques are used to predict the failure
life more precisely using steady-state solution (\citet{Leckie1974},
\citet{Nikitenko1979}, \citet{Hyde1998,Hyde1999}, \citet{Perrin2000}).
As shown in \citet{Browne1981}, \citet{Coleman1985} for the particular pipe
welds investigated, the steady-state analysis underestimates the
failure time by about 20-40 percent, but predicts failure position
quite well.

To simplify the parametric analysis, we study a family of
boundary-value problems for multi-component pipe creep depending
on a small parameter $s$. When $s=0$, the problem is reduced to
the case of homogeneous pipe creep. The corresponding solution
$\boldsymbol{\sigma}(0)$ for the steady-state stress distribution
is well known (\citet{Odqvist1974}, \citet{Malinin1981}, \citet{Hyde1996}) and
considered to be the basic solution. In order to get the common
solution $\boldsymbol \sigma(s)$ from this basic one,
correction terms should be added.
The equations for the correction terms are formally obtained as
perturbations of a boundary-value problem with respect to $s$.
These equations formulate a problem of linear elasticity with
respect to linear-elastic solid with anisotropic elastic
properties.

The utility of used technique is guaranteed especially because the
theory of linear elasticity is in a very satisfactory state of
completion; every complicated case of parameter distribution can
be treated in a routine manner as a combination of simple ones.

Correction terms are obtained numerically for some problems with
the help of the Ritz method. In some cases the simplicity of
geometry enables us to construct an approximate analytical
solution. An exact expression for stress jumps at the interface is obtained.
Numerical solutions of the nonlinear problem are
obtained with the help of ANSYS finite element code for
comparison.

Throughout light-face letters we denote scalars, the
bald-face letters stand for tensors.
The notation $\vec{( \ )}$ is used in vector-matrix
form of constitutive equations to designate vectors
$\vec{\boldsymbol \sigma}$ and  $\vec{\boldsymbol \varepsilon}$
of stress and strain components.

\section{Problem formulation: two-material model of narrow gap weldment}
In this section we consider a two-material model only. It will be
shown later that the solution for some multi-material models can be
reduced to this case.

\begin{figure}
\psfrag{A}[m][][1][0]{parent} \psfrag{AA}[m][][1][0]{material}
\psfrag{C}[m][][1][0]{$H$} \psfrag{B}[m][][1][0]{weldment}
\psfrag{D}[m][][1][0]{$h$} \psfrag{E}[m][][1][0]{(a)}
\psfrag{F}[m][][1][0]{(b)} \psfrag{G}[m][][1][0]{$r_o$}
\psfrag{H}[m][][1][0]{$r_i$} \psfrag{K}[m][][1][0]{$\theta$}
\psfrag{L}[m][][1][0]{$r$} \psfrag{M}[m][][1][0]{$z$}
\psfrag{N}[m][][1][0]{$\sigma_{rz}=0$}
\psfrag{O}[m][][1][0]{$u_z=0$}
\psfrag{P}[m][][1][0]{$\sigma_{r}=-p$}
\psfrag{Q}[m][][1][0]{$\sigma_{rz}=0$}
\psfrag{R}[m][][1][0]{$\sigma_{r}=0$}
\scalebox{1}{\includegraphics{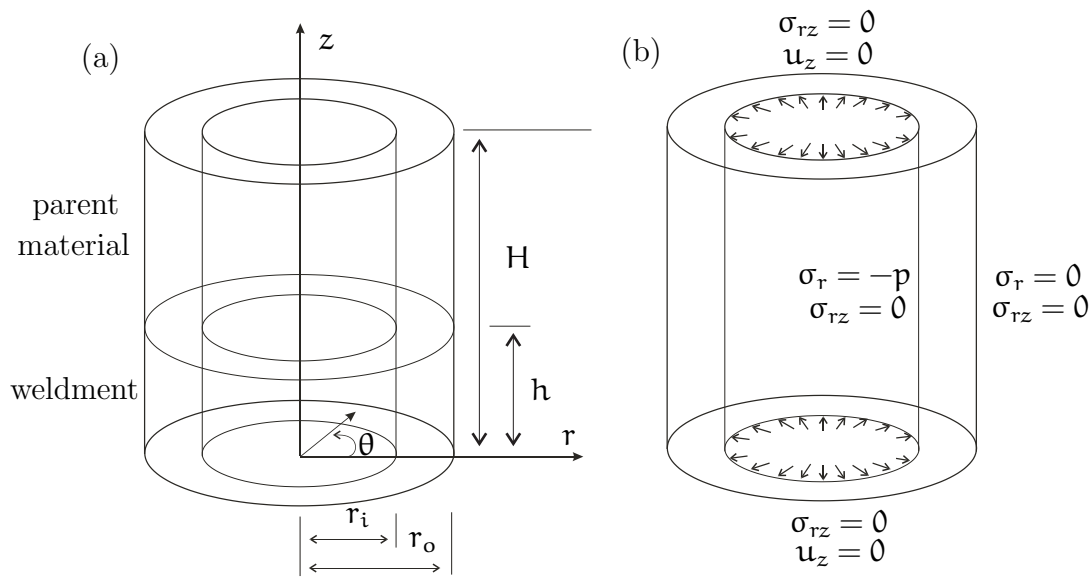}} \caption{System
configuration (a) and boundary conditions (b) \label{fig1}}
\end{figure}

The configuration analyzed is shown in Fig. \ref{fig1}.
Assume that $(r,z,\theta)\in\Omega\times[0,2\pi)$,
$(r,z,\theta)\in\Omega^{-}\times[0,2\pi)$,
$(r,z,\theta)\in\Omega^{+}\times[0,2\pi)$ describe the volume
occupied by the solid, by the weld metal, and by the parent
material, respectively. Here $\Omega^{-}=[r_i,r_o]\times[0,h)$,
$\Omega^{+}=[r_i,r_o]\times[h,H]$,
$\Omega=\Omega^{-}\cup\Omega^{+}$. The basic equations of the
problem are given below.

1. Equations of equilibrium

\begin{equation}\label{equilInvar}
 \boldsymbol \nabla \cdot \boldsymbol \sigma = \boldsymbol 0,
\end{equation}
where $\boldsymbol \sigma$ is the stress tensor. In
\eqref{equilInvar} the volumetric forces are ignored.

 2. Strain-displacement relations
\begin{equation}\label{str-dispInvar}
\boldsymbol \varepsilon =\frac{1}{2}  \left(  \boldsymbol
\nabla \boldsymbol u + \boldsymbol \nabla
\boldsymbol u^T \right),
\end{equation}
where $\boldsymbol \varepsilon$ is the linearized strain tensor
and $\boldsymbol u$ is the displacement vector.

3. The governing equations can be summarized as Norton's creep law
(see, for example, \citet{Odqvist1974})
\begin{equation}\label{norton}
\dot{\boldsymbol\varepsilon} = A \ \boldsymbol s \ (\sigma_{vM})^{n-1},
\quad \boldsymbol s = \boldsymbol \sigma - \frac{1}{3}tr(\boldsymbol \sigma) \boldsymbol I, \quad
\sigma_{vM}=\sqrt{\frac{3}{2}\boldsymbol s : \boldsymbol s},
\end{equation}
where $\dot{( \ )}$ is the time derivative, $\boldsymbol s$ is the stress
deviator, $\sigma_{vM}$ is the von Mises equivalent stress, $\boldsymbol
I$ is the second rank unit tensor, and $A$, $n$ are material
constants. Note that Norton's creep law
is often written as
$\dot{\boldsymbol\varepsilon} = \frac{3}{2} a \ \boldsymbol s \ (\sigma_{vM})^{n-1}$.
It means $A= \frac{3}{2} a$ in our notation.

If we take into account that the problem is axisymmetric, we have
in cylindrical coordinates the following equilibrium equations

\begin{equation}\label{equil}
 \frac{\partial \sigma_{r}}{\partial r}+\frac{1}{r}(\sigma_{r}-\sigma_{\theta})+
 \frac{\partial \sigma_{rz}}{\partial z}=0, \  \
\frac{\partial \sigma_{rz}}{\partial r}+\frac{1}{r}\sigma_{rz}+
 \frac{\partial \sigma_{z}}{\partial z}=0,
\end{equation}

strain-displacement relations

\begin{equation}\label{str-disp}
\varepsilon_{r}=\frac{\partial u_r}{\partial r}, \
\varepsilon_{\theta}=\frac{u_r}{r}, \
 \varepsilon_{z}=\frac{\partial u_z}{\partial z}, \
 \varepsilon_{rz}=\frac{1}{2}  \left( \frac{\partial u_r}{\partial z}+\frac{\partial u_z}{\partial r} \right),
\end{equation}

and boundary conditions
\begin{equation}
\sigma_{r}=-p, \ \sigma_{rz}=0 \quad \text{at} \quad r=r_i; \
\sigma_{r}=\sigma_{rz}=0  \quad \text{at} \quad r=r_o;
\end{equation}
\begin{equation*}
\sigma_{rz}=0, \ u_z=0 \quad \text{at} \quad z=0, \ z=H.
\end{equation*}

Finally we define the distribution of the parameter $A(r,z)$ in
Norton's creep law \eqref{norton} as a piecewise constant function in
$\Omega$
\begin{equation}\label{property0}
A=A^+ \ \ \text{in} \ \Omega^{+}, \ A=A^- \ \ \text{in} \
\Omega^{-}.
\end{equation}

The stress field $\boldsymbol \sigma$, defined by
\eqref{equilInvar} --- \eqref{property0}, does not change if
instead of $A(r,z)$ we use $\lambda A(r,z)$, where $\lambda$ is
any nonzero constant. Consequently, without loss of generality it
can be assumed that\footnote{Here $s$ is not necessary a small
parameter. In practice, $A^+$ and $A^-$ might differ essentially.}
\begin{equation}\label{property}
A=A^+=1 \ \ \text{in} \ \Omega^{+}, \ A=A^-=1-s \ \ \text{in} \
\Omega^{-}.
\end{equation}

If we eliminate displacements from the strain-displacement
relations \eqref{str-disp}, we obtain
compatibility equations
\begin{equation}\label{compat}
C_1(\boldsymbol\varepsilon)=C_2(\boldsymbol\varepsilon)=0,
\end{equation}
\begin{equation}\label{compat2}
C_1(\boldsymbol\varepsilon)=r(\varepsilon_r-\frac{\partial(r
\varepsilon_\theta)}{\partial r}),  \ C_2(\boldsymbol\varepsilon)=r
\frac{\partial^2 \varepsilon_\theta}{\partial z^2}+
\frac{\partial\varepsilon_z}{\partial
r}-2\frac{\partial\varepsilon_{rz}}{\partial z}.
\end{equation}

We consider the weak form of compatibility equations
expressed by the equation of complementary virtual power principle
\citep{Washizu1982}
\begin{equation}\label{comvpow}
L(\boldsymbol\sigma,s)<\delta\boldsymbol\sigma>=0, \quad \forall \delta\boldsymbol\sigma,
\end{equation}
\begin{equation*}
L(\boldsymbol\sigma,s)<\delta\boldsymbol\sigma>\equiv\int\limits_{\Omega}\dot{\boldsymbol\varepsilon}(\boldsymbol\sigma,s)
: \delta\boldsymbol\sigma d\Omega.
\end{equation*}
Here we use the brackets $< \cdot >$ to enclose the argument of
a linear operator; $\dot{\boldsymbol\varepsilon}(\boldsymbol\sigma,s)$ is the
strain rate defined by \eqref{norton} and \eqref{property};
$d\Omega=r dr dz$; $\delta\boldsymbol\sigma$ is a virtual stress field
that satisfy the equations of equilibrium \eqref{equil} and
homogeneous boundary conditions
\begin{equation}\label{homog}
\sigma_{r}=\sigma_{rz}=0 \quad \text{at} \quad r=r_i, \ r=r_o; \
\sigma_{rz}=0 \quad  \text{at} \quad z=0, \ z=H.
\end{equation}

In what follows we search function $\boldsymbol\sigma(s)$, such that
\begin{equation}\label{impcur}
L(\boldsymbol\sigma(s),s)<\delta\boldsymbol\sigma>=0, \quad \forall
\delta\boldsymbol\sigma.
\end{equation}

\section{Perturbation method}
\subsection{Unperturbed solution: creep response of a homogeneous pipe}
Now suppose that $A^+=A^-$, i.e. $s=0$. The problem
\eqref{comvpow} is reduced to a one-dimensional, and the solution
$\boldsymbol\sigma^0$ of this problem is well known (\citet{Odqvist1974}, \citet{Malinin1981}, \citet{Hyde1996})
\begin{equation}\label{basic}
\sigma^0_{r}=a+a_r r^{-2/n}, \ \sigma^0_{\theta}=a+a_\theta
r^{-2/n}, \ \sigma^0_{z}=a+a_z r^{-2/n},  \ \sigma^0_{rz}=0,
\end{equation}
\begin{equation*}
a=p \ \frac{r^{2/n}_i}{r^{2/n}_o-r^{2/n}_i}, \ a_r= - p \ \frac{r^{2/n}_i \
r^{2/n}_o}{r^{2/n}_o-r^{2/n}_i}, \ a_\theta= \frac{n-2}{n} \ a_r,
\ \ a_z= \frac{n-1}{n} \ a_r.
\end{equation*}

\subsection{First perturbation}
Let $s$ be a small parameter. Assume that\footnote{The
justification of perturbation method could be performed by means
of Implicit Function Theorem \citep{Antman1995}. Precise conditions
justifying perturbation method should be formulated in an
appropriate function space.}

\begin{equation}\label{hypotes}
\boldsymbol\sigma(s)=\boldsymbol\sigma^0+ s \boldsymbol\sigma^1 +o(s).
\end{equation}
Here $\boldsymbol\sigma^1=\frac{d\boldsymbol\sigma(s)}{ds}\mid_{s=0}$
is the unknown derivative which must satisfy the equations of
equilibrium \eqref{equil} and the homogeneous boundary conditions
\eqref{homog}; $o(s)$ is the little-o Landau symbol. From
\eqref{property} it follows that
\begin{equation}\label{decomp}
L(\boldsymbol\sigma,s)<\delta\boldsymbol\sigma>=\left(L^0(\boldsymbol\sigma)+ s
L^1(\boldsymbol\sigma)\right)<\delta\boldsymbol\sigma>
\end{equation}
with $L^0(\boldsymbol\sigma)\equiv L(\boldsymbol\sigma,0)$,
$L^1(\boldsymbol\sigma)\equiv\frac{\partial L(\boldsymbol\sigma,s)}{\partial s}$.
Substituting \eqref{hypotes} and \eqref{decomp}
in \eqref{impcur}, we get
\begin{equation*}
\left(L^0(\boldsymbol\sigma^0)+s
\frac{dL^0}{d\boldsymbol\sigma}\mid_{\boldsymbol\sigma=\boldsymbol\sigma^0}<\boldsymbol\sigma^1>+
s L^1(\boldsymbol\sigma^0)+ o(s)\right)<\delta\boldsymbol\sigma>=0.
\end{equation*}
Since $L^0(\boldsymbol\sigma^0)<\delta\boldsymbol\sigma>=0$ and $\frac{o(s)}{s}
<\delta\boldsymbol\sigma> \rightarrow 0$ as $s\rightarrow 0$, we have an
equation for $\boldsymbol\sigma^1$
\begin{equation}\label{firper}
\frac{dL^0}{d\boldsymbol\sigma}\mid_{\boldsymbol\sigma=\boldsymbol\sigma^0}<\boldsymbol\sigma^1><\delta\boldsymbol\sigma>=-L^1(\boldsymbol\sigma^0)<\delta\boldsymbol\sigma>,
\ \forall \delta\boldsymbol\sigma.
\end{equation}
We will analyze more closely the linear operator
$\frac{dL^0}{d\boldsymbol\sigma}\mid_{\boldsymbol\sigma=\boldsymbol\sigma^0}<\cdot>$ in
the next section.

\section{Auxiliary problem}
\subsection{Linear elastic material}
In the previous section it was shown that the correction term
$\boldsymbol\sigma^1$ can be found from the linear equation
\eqref{firper}. It is clear that
\begin{equation}\label{derivchange}
\frac{dL^0}{d\boldsymbol\sigma}\mid_{\boldsymbol\sigma=\boldsymbol\sigma^0}<\boldsymbol
\sigma^1><\delta\boldsymbol\sigma>=
\int\limits_{\Omega}\frac{\partial\dot{\boldsymbol\varepsilon}}{\partial
\boldsymbol\sigma}(\boldsymbol\sigma^0,0)<\sigma^1> : \delta\boldsymbol\sigma d\Omega,
\end{equation}
where
\begin{equation}\label{eps}
\dot{\boldsymbol\varepsilon}(\boldsymbol\sigma,0)=\boldsymbol s(\boldsymbol\sigma) \
(\sigma_{vM}(\boldsymbol\sigma))^{n-1}.
\end{equation}
Let us introduce a vector notation $\vec{( \ )}$ as
\begin{equation}\label{notation}
\vec{\boldsymbol \sigma}=(\sigma_r, \sigma_{\theta}, \sigma_z,
\sigma_{rz})^T, \ \vec{\boldsymbol \varepsilon}=(\varepsilon_r,
\varepsilon_{\theta}, \varepsilon_z, 2\varepsilon_{rz})^T.
\end{equation}
Substituting \eqref{notation} in \eqref{derivchange} and
differentiating \eqref{eps} we obtain
\begin{equation}\label{tezor}
\frac{dL^0}{d\boldsymbol\sigma}\mid_{\boldsymbol\sigma=\boldsymbol\sigma^0}<\boldsymbol
\sigma^1><\delta\boldsymbol\sigma>= \int\limits_{\Omega} (\vec{\boldsymbol
\sigma^1})^T \ \boldsymbol C \ \vec{\delta \boldsymbol\sigma} \ d\Omega.
\end{equation}
Here we have introduced the compliance matrix as follows
\arraycolsep10pt
\begin{equation}\label{tezorC}
\boldsymbol C=\frac{\partial \vec{\dot{\boldsymbol \varepsilon}}}{\partial
\vec{\boldsymbol \sigma}}(\vec{\boldsymbol \sigma^0})=
\end{equation}
\begin{equation*}
\left(\sqrt{3} \frac{|a_r|}{n}r^{  \displaystyle  -2/n}\right)^{
\displaystyle  n-2} \left[\begin{array}{cccc}
       2/3 + (n-1)/2 & -1/3-(n-1)/2 &  -1/3 & 0 \\
     &  2/3 + (n-1)/2      &  -1/3 & 0 \\
     sym.  &  &    -1/3  &     0  \\
          &     &      &     2  \\
\end{array} \right] \,.
\end{equation*}

Thus, the left-hand side of \eqref{firper} can be treated as an
internal complementary virtual work
with respect to linear elastic solid with constitutive law
\begin{equation}\label{linlaw}
(\varepsilon_r, \varepsilon_{\theta}, \varepsilon_z,
2\varepsilon_{rz})^T=\boldsymbol C (\sigma_r, \sigma_{\theta}, \sigma_z,
\sigma_{rz})^T.
\end{equation}

The eigenvalues of the compliance operator $\boldsymbol C$ are
\begin{equation}\label{tezor}
(\lambda_1,\lambda_2,\lambda_3,\lambda_4)=(\sqrt{3}
\frac{|a_r|}{n}r^{-2/n})^{n-2} (0,2,2,n).
\end{equation}
The problem of steady-state creep is reduced to the elasticity
problem for orthotropic, incompressible, and inhomogeneous solid
\eqref{linlaw}.

\subsection{Displacement jump}
We now seek to convert the right-hand side of \eqref{firper} into
a surface integral through Gauss theorem. It can be proved that
\begin{equation}\label{load}
-L^1(\boldsymbol\sigma^0)<\delta\boldsymbol\sigma>=
\int\limits_{\Omega^{-}}\dot{\boldsymbol\varepsilon}(\boldsymbol\sigma^0,0)
: \delta\boldsymbol\sigma \ d \Omega=
 - \int\limits_{\{z=h\}} 3^{\frac{n-1}{2}}a_r |a_r|^{n-1}\frac{1}{r n^n}\delta\sigma_{rz} r dr.
\end{equation}

\begin{figure}
\psfrag{A}[m][][1][0]{$0$} \psfrag{B}[m][][1][0]{$z$}
\psfrag{C}[m][][1][0]{$r$} \psfrag{D}[m][][1][0]{$u^+_r-u^-_r$}
\psfrag{E}[m][][1][0]{$\Omega^-$}
\psfrag{F}[m][][1][0]{$\Omega^+$} \psfrag{G}[m][][1][0]{$H$}
\psfrag{H}[m][][1][0]{$h$}
\scalebox{1}{\includegraphics{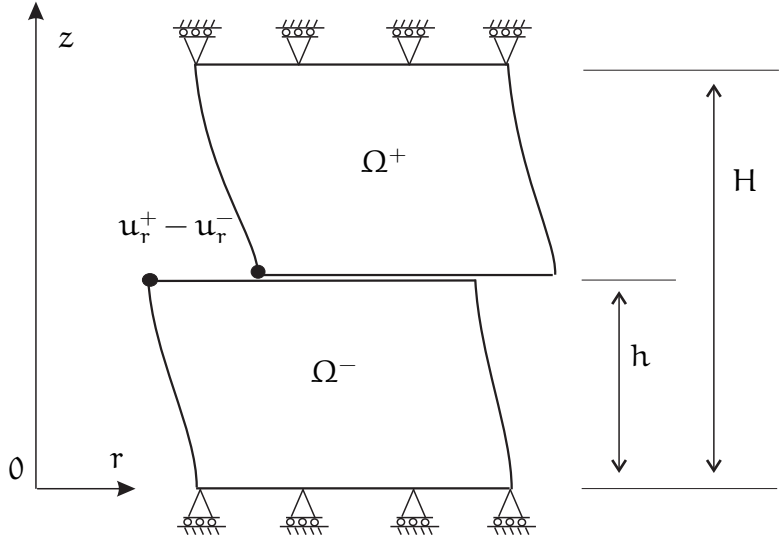}} \caption{Displacement
jump \label{fig2}}
\end{figure}

Let us show that \eqref{load} prescribes a jump of displacements
at the interface $\{z=h\}$ (see fig. \ref{fig2}).
Consider the principle of virtual complementary work
for linear elastic solids $\Omega^{-}$ and $\Omega^{+}$. One gets

\begin{equation}\label{twosolids1}
\int\limits_{\Omega^-} (\vec{\boldsymbol \sigma^-})^T \ \boldsymbol C \
\vec{\delta \boldsymbol\sigma} \ d\Omega=
 \int\limits_{\{z=h\}} (u^-_r\delta\sigma_{rz}+u^-_z\delta\sigma_{z}) r dr,
\end{equation}
\begin{equation}\label{bound1}
\delta\sigma_{r}=\delta\sigma_{rz}=0 \ \text{at} \ r=r_i, \ r=r_o;
\ \delta\sigma_{rz}=0 \ \text{at} \ z=0.
\end{equation}
\begin{equation}\label{twosolids2}
\int\limits_{\Omega^+} (\vec{\boldsymbol \sigma^+})^T \ \boldsymbol C \
\vec{\delta \boldsymbol\sigma} \ d\Omega=
 -\int\limits_{\{z=h\}} (u^+_r\delta\sigma_{rz}+u^+_z\delta\sigma_{z}) r dr,
\end{equation}
\begin{equation}\label{bound2}
\delta\sigma_{r}=\delta\sigma_{rz}=0 \ \text{at} \ r=r_i, \ r=r_o;
\ \delta\sigma_{rz}=0 \ \text{at} \ z=H.
\end{equation}
Suppose
\begin{equation}\label{stresscont}
\sigma^+_{rz}-\sigma^-_{rz}=0, \quad \sigma^+_{z}-\sigma^-_{z}=0
\ \text{at}  \ z=h;
\end{equation}
\begin{equation}\label{dispjump}
u^+_z-u^-_z=0, \quad u^+_r-u^-_r=3^{\frac{n-1}{2}}a_r
|a_r|^{n-1}\frac{1}{r n^n} \ \text{at}  \ z=h;
\end{equation}
then it can be shown that the solution $\boldsymbol\sigma^1$ of
\eqref{firper} has the form
\begin{equation}\label{decomposition}
\boldsymbol\sigma^1=\begin{cases}
    \boldsymbol\sigma^+ \quad \text{in} \ \Omega^+\\
    \boldsymbol\sigma^- \quad \text{in} \ \Omega^-
  \end{cases}.
\end{equation}

\subsection{Stress jump}
We define a notation for jumps of field variables at the interface
${\{z=h\}}$
\begin{equation*}
[v]=v^+ - v^-.
\end{equation*}
Substituting \eqref{dispjump} for $u_r$ in \eqref{str-disp}, we
get
\begin{equation*}
[\varepsilon_r]=-c/r^2, \quad [\varepsilon_\theta]=c/r^2, \quad
c=3^{\frac{n-1}{2}}a_r |a_r|^{n-1}\frac{1}{n^n}.
\end{equation*}
If we combine this with \eqref{tezorC}, \eqref{linlaw}, and
\eqref{stresscont}, we obtain
\begin{equation}\label{strjump}
[\sigma^1_r]=-\frac{\sqrt{3} a_r |a_r|}{n^3} r^{-4/n}, \quad
[\sigma^1_{\theta}]=\frac{\sqrt{3} a_r |a_r|}{n^3} r^{-4/n}.
\end{equation}

\subsection{Case of multi-material structure}
\begin{figure}
\psfrag{A}[m][][1][0]{$z_1=0$} \psfrag{B}[m][][1][0]{$z_2$}
\psfrag{C}[m][][1][0]{$z_3$} \psfrag{D}[m][][1][0]{$z_4=H$}
\psfrag{E}[m][][1][0]{$z$}
\psfrag{F}[m][][1][0]{$A_1=1-s(1+\alpha_1)$}
\psfrag{G}[m][][1][0]{$A_2=1-s$} \psfrag{H}[m][][1][0]{$A_3=1$}
\scalebox{1}{\includegraphics{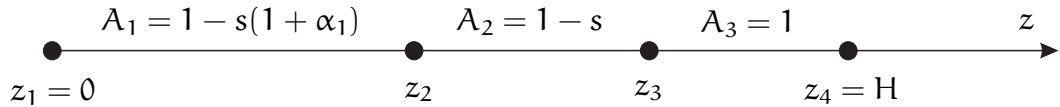}} \caption{Distribution of
$A$ in $m$-material structure ($m=3$) \label{fig3}}
\end{figure}

In this subsection we consider $m$-material structure ($m \geq 2$)
\begin{equation*}
\Omega = \Omega_1 \cup \Omega_2 \cup ... \cup \Omega_m,
\end{equation*}
\begin{equation*}
\Omega_j=[r_i,r_o] \times [z_j, z_{j+1}], \ z_1=0, \ z_{m+1}=H,
\end{equation*}
\begin{equation*}
A=A_j \ \text{in} \ \Omega_j, \ j=1,...,m; \ A_m=1, \ A_{m-1} \neq
A_m.
\end{equation*}
There exist unique $s$ and $\alpha_j, j=1,...,m-1$ such that (fig.
\ref{fig3})
\begin{equation}\label{sum}
A(z)=1-s \ \sum^{m-1}_{j=1} \alpha_j \ \text{H}(z_{j+1}-z), \
\alpha_{m-1}=1.
\end{equation}
Here $\text{H}(z)$ is the Heaviside function.

If we use the small parameter method to solve this problem, then
we get an approximation in the form \eqref{hypotes}; $\boldsymbol\sigma^0$
is given by \eqref{basic} and $\boldsymbol\sigma^1$ is obtained from
\eqref{firper}. Since \eqref{sum} is valid, it follows that the
right-hand side of \eqref{firper} has the form
\begin{equation*}
-L^1(\boldsymbol\sigma^0)<\delta\boldsymbol\sigma>= \sum^{m-1}_{j=1} \alpha_j \
 \int\limits_{\{z=z_j\}} \frac{-c}{r}\delta\sigma_{rz} r dr.
\end{equation*}
Owing to the fact that \eqref{firper} is linear, the perturbation
$\boldsymbol\sigma^1$ is a linear combination of $m-1$ solutions of type
\eqref{decomposition}.

\section{Solution techniques}

\subsection{Approximate analytical solution}

In this section we construct an analytical solution for the first
perturbation term $\boldsymbol \sigma^1$. Consider a pair of stress
functions $(\varphi,\psi)$, such that condition \eqref{equil} is
satisfied
\begin{equation}\label{strfunc}
\sigma_r=\varphi, \ \sigma_{\theta}=\frac{\partial (r
\varphi)}{\partial r}+\frac{\partial^2 \psi}{\partial z^2}, \
\sigma_z=-\frac{1}{r}\frac{\partial \psi}{\partial r},  \
\sigma_{rz}=\frac{1}{r}\frac{\partial \psi}{\partial z}.
\end{equation}
These functions were constructed from compatibility equations
\eqref{compat}, \eqref{compat2} using the standard technique
\citep{Washizu1982}. This choice of stress functions simplifies the
application of variational methods.

We use the Kantorovich method (see, for example,
\cite{Kantorovich1958}) to reduce the 2-D variational problem to 1-D
variational problem. Suppose that for $(\varphi,\psi)$, that
define the solution $\boldsymbol \sigma^1$ of \eqref{firper}, the
following static hypothesis is valid
\begin{equation}\label{hypothesis}
\varphi(r,z)=\varphi_1(r) \varphi_2(z), \ \psi(r,z)=\psi_1(r)
\psi_2(z),
\end{equation}
where $\varphi_1(r), \psi_1(r)$ are given Kantorovich trial
functions, such that
\begin{equation*}
\varphi_1(r_i)=\varphi_1(r_o)=\psi_1(r_i)=\psi_1(r_o)=0.
\end{equation*}
We use the Kantorovich method  to obtain a system of differential
equations and boundary conditions for $\varphi_2(z), \psi_2(z)$.
This method gives a projection of $\boldsymbol \sigma^1$ on the subspace,
defined by \eqref{hypothesis}. It is more convenient to solve
\eqref{firper} in the form \eqref{twosolids1} ---
\eqref{decomposition}. For the sake of brevity we describe only
the application of the Kantorovich method for equation
\eqref{twosolids1}. For equation \eqref{twosolids2} this procedure
can be arranged in a similar manner.

Since $\varphi_1(r)$ and $\psi_1(r)$ are fixed, the variation of
\eqref{hypothesis} gives
\begin{equation}\label{hypothesis2}
\delta \varphi=\varphi_1(r) \delta \varphi_2(z), \ \delta
\psi=\psi_1(r) \delta \psi_2(z).
\end{equation}
Let us rewrite \eqref{twosolids1}, \eqref{bound1} in terms of
$\varphi_2,\psi_2$
\begin{equation*}
\int\limits_{\Omega^-} (\vec{\boldsymbol \sigma}(\varphi_2,\psi_2))^T \
\boldsymbol C \ \vec{\boldsymbol \sigma}(\delta\varphi_2,\delta\psi_2) \ d\Omega=
\end{equation*}
\begin{equation*}
 \int^{r_o}_{r_i} (u^-_r \ \frac{\psi_1(r)}{r} \frac{d \delta\psi_2}{dz}(h)-u^-_z\
 \frac{\delta \psi_2(h)}{r} \frac{d \psi_1(r)}{dr}) r dr,
\end{equation*}
\begin{equation*}
\frac{d \psi_2}{dz}(0)=0, \quad \frac{d \delta \psi_2}{dz}(0)=0.
\end{equation*}
Here $\vec{\boldsymbol \sigma}(\varphi_2,\psi_2)$ and $\vec{\boldsymbol
\sigma}(\delta\varphi_2,\delta\psi_2)$ are defined by
\eqref{notation}, \eqref{hypothesis}, and \eqref{hypothesis2}.
Using Gauss theorem and the fundamental lemma of the calculus of
variation \citep{Washizu1982} we obtain compatibility equations and
boundary conditions
\begin{equation}\label{eq1}
\int^{r_o}_{r_i} \varphi_1(r) C_1(\vec{\boldsymbol
\varepsilon}(\varphi_2,\psi_2)) dr=0, \quad \int^{r_o}_{r_i}
\psi_1(r) C_2(\vec{\boldsymbol \varepsilon}(\varphi_2,\psi_2)) dr=0,
\end{equation}
\begin{equation}\label{eq2}
\int\limits_{\{z=0\}}(2 \varepsilon_{rz}-\frac{\partial (r
\varepsilon_{\theta})}{\partial z}) \psi_1(r) dr=0,
\int\limits_{\{z=h\}}(2 \varepsilon_{rz}-\frac{\partial (r
\varepsilon_{\theta})}{\partial z}- \frac{d u^-_z}{d r})
\psi_1(r) dr=0,
\end{equation}
\begin{equation}\label{eq22}
\int\limits_{\{z=h\}}(\varepsilon_{\theta}- \frac{u^-_r}{r})
\psi_1(r) dr=0.
\end{equation}
Here $\vec{\boldsymbol \varepsilon} = \boldsymbol C \ \vec{\boldsymbol \sigma}$ and
$C_1(\vec{\boldsymbol \varepsilon}), C_2(\vec{\boldsymbol \varepsilon})$ are
defined by \eqref{compat2}, \eqref{notation}.  After integration,
system \eqref{eq1} has the form
\begin{equation}\label{eq3}
a_1 \varphi_2+a_2 \psi_2+a_3 \frac{d^2 \psi_2}{d z^2}=0, \ b_1
\varphi_2 + b_2 \frac{d^2 \varphi_2}{d z^2} + b_3 \psi_2 + b_4
\frac{d^2 \psi_2}{d z^2} + b_5 \frac{d^4 \psi_2}{d z^4}=0.
\end{equation}
Elimination of $\varphi_2$ from \eqref{eq3} gives
\begin{equation}\label{eq4}
k_1 \psi_2 + k_2 \frac{d^2 \psi_2}{d z^2} + k_3 \frac{d^4
\psi_2}{d z^4}=0.
\end{equation}
After integration in \eqref{eq2}, \eqref{eq22}, we obtain
\begin{equation*}
(e_1 \frac{d \psi_2}{d z}+ e_2 \frac{d^2 \psi_2}{d z^2} + e_3
\frac{d^3 \psi_2}{d z^3})(0)=0,
\end{equation*}
\begin{equation*}
(e_1 \frac{d \psi_2}{d z}+ e_2 \frac{d^2 \psi_2}{d z^2} + e_3
\frac{d^3 \psi_2}{d z^3})(h)= \int\limits_{\{z=h\}} \frac{d
u^-_z}{d r} \psi_1(r) dr,
\end{equation*}
\begin{equation*}
(g_1 \psi_2 + g_2 \frac{d^2 \psi_2}{d
z^2})(h)=\int\limits_{\{z=h\}} \frac{u^-_r}{r} \psi_1(r) dr.
\end{equation*}
Dealing with \eqref{twosolids2}, \eqref{bound2} in a similar
fashion we obtain compatibility equations and boundary conditions
for $\psi_2(z)$ in $(h,H)$. Finally, taking into account
\eqref{stresscont}, \eqref{dispjump}, we have eight boundary
conditions
\begin{equation*}
\frac{d \psi_2}{dz}(0)=0, \quad \frac{d \psi_2}{dz}(H)=0,
\end{equation*}
\begin{equation*}
[\psi_2]=0, \quad [\frac{d \psi_2}{dz}]=0,
\end{equation*}
\begin{equation*}
(e_1 \frac{d \psi_2}{d z}+ e_2 \frac{d^2 \psi_2}{d z^2} + e_3
\frac{d^3 \psi_2}{d z^3})(0)=0,
\end{equation*}
\begin{equation*}
(e_1 \frac{d \psi_2}{d z}+ e_2 \frac{d^2 \psi_2}{d z^2} + e_3
\frac{d^3 \psi_2}{d z^3})(H)=0,
\end{equation*}
\begin{equation*}
[e_1 \frac{d \psi_2}{d z}+ e_2 \frac{d^2 \psi_2}{d z^2} + e_3
\frac{d^3 \psi_2}{d z^3}]=0,
\end{equation*}
\begin{equation*}
[g_1 \psi_2 + g_2 \frac{d^2 \psi_2}{d z^2}]=\int^{r_o}_{r_i}
\frac{c}{r^2} \psi_1(r) dr.
\end{equation*}
Solution $\psi_2(z)$ is uniquely defined by these equations
combined with two ordinary differential equations (ODE)
 \eqref{eq4} of fourth order (one equation in $(0,h)$ and another in $(h,H)$).

\subsection{Numerical solution technique}
In this section we solve problem \eqref{firper} numerically with
the help of the Ritz method. Let $\boldsymbol \sigma_i, i=1,2,...,N$ be a
system of stresses, satisfying \eqref{equil} and \eqref{homog}.
Suppose
\begin{equation*}
\boldsymbol \sigma^1=\sum^{N}_{i=1} c_i \boldsymbol \sigma_i.
\end{equation*}
Then the unknown constants $c_i$ are defined from a system of $N$
linear algebraic equations
\begin{equation*}
\frac{dL^0}{d\boldsymbol\sigma}\mid_{\boldsymbol\sigma=\boldsymbol\sigma^0}<\boldsymbol\sigma^1><\boldsymbol
\sigma_i>=-L^1(\boldsymbol\sigma^0)<\boldsymbol \sigma_i>, \ i=1,...,N.
\end{equation*}
We use the system $\boldsymbol \sigma_i, i=1,2,...$, produced by means of
\eqref{strfunc}. The complete system of $\varphi, \psi$ is given
by combinations of trigonometric functions
\begin{equation}\label{ser1}
\varphi=\sin\left(i \frac{r-r_i}{r_o-r_i} \pi\right) \ \cos\left(j \frac{ \ z \
\pi}{H}\right), \ i=1,2,...,N_r , \ j=0,1,...,N_z ,
\end{equation}
\begin{equation}
\varphi= \sin\left(i \frac{r-r_i}{r_o-r_i} \pi\right) \ \cos\left(\frac{z \ \pi}{2
H}\right), \ i=1,2,...,N_r ,
\end{equation}
\begin{equation}
\varphi= \sin\left(i \frac{r-r_i}{r_o-r_i} \pi\right) \ \sin\left(\frac{z \ \pi}{2
H}\right), \ i=1,2,...,N_r ,
\end{equation}
\begin{equation}
\psi= \sin\left(n \frac{r-r_i}{r_o-r_i} \pi\right) \ \cos\left(k \frac{ \ z \
\pi}{H}\right),\ n=1,2,...,N_r , \ k=0,1,...,N_z.
\end{equation}

In order to approximate the stress-jump \eqref{strjump} at the
interface $\{z=h\}$ it could be useful to consider additionally
discontinuous functions
\begin{equation}
\varphi=\sin\left(i \frac{r-r_i}{r_o-r_i} \pi\right) \ \text{H}(z), \
i=1,2,...,N_r,
\end{equation}
\begin{equation}\label{ser2}
\psi= \sin\left(n \frac{r-r_i}{r_o-r_i} \pi\right)  \ \begin{cases}
  \displaystyle  \frac{z^2}{2h}, \ z < h\\
   \displaystyle \frac{-z^2}{2(H-h)}+ \frac{z H}{H-h}+\frac{-h H}{2(H-h)}, \  z \geq h
  \end{cases},  n=1,2,...,N_r.
\end{equation}
Here $\text{H}(z)$ is the Heaviside function.

If we put $N_r=1, \ N_z \rightarrow \infty$, then the solution
tends to the approximate analytical solution from the previous
section.

\section{Results of calculations}
The specimen dimensions, applied loads, and material constants are
given by $H=8$, $h=0.5$, $r_i=1$, $r_o=2$, $p=1$, $n=3$.

\subsection{Solution of the linear auxiliary problem}
First we compare the approximate analytical solution of
\eqref{firper} (section 5.1) with the numerical solution which
was obtained using the Ritz method (section 5.2).

We used $\varphi_1(r)=\psi_1(r)=\sin(\frac{r-r_i}{r_o-r_i} \pi)$
in \eqref{hypothesis} to construct the analytical solution. Values
of constants $a_i, b_i, k_i, e_i, g_i$ are given in Appendix. The
numerical solution by the Ritz method is performed witn
 $N_r=N_z=25$ in \eqref{ser1} --- \eqref{ser2}.

\begin{figure}
\psfrag{A}[m][][1][0]{$0$} \psfrag{B}[m][][1][0]{$1$}
\psfrag{C}[m][][1][0]{$2$} \psfrag{D}[m][][1][0]{$3$}
\psfrag{E}[m][][1][0]{$4$} \psfrag{F}[m][][1][0]{$5$}
\psfrag{G}[m][][1][0]{$6$} \psfrag{H}[m][][1][0]{$7$}
\psfrag{K}[m][][1][0]{$8$} \psfrag{L}[m][][1][0]{$z$}
\psfrag{M}[m][][1][0]{$\sigma_{rz}$}
\psfrag{N}[m][][1][0]{$-0.05$} \psfrag{O}[m][][1][0]{$0.05$}
\psfrag{P}[m][][1][0]{$0.1$} \psfrag{R}[m][][1][0]{$\sigma_{vM}$}
\psfrag{Q}[m][][1][0]{$r=1.5$} \psfrag{S}[m][][1][0]{$r=1.2$}
\psfrag{T}[m][][1][0]{$0.2$} \psfrag{U}[m][][1][0]{$-0.1$}
\psfrag{X}[m][][1][0]{$\sigma_{\theta}$}
\psfrag{Y}[m][][1][0]{$-0.2$} \psfrag{Z}[m][][1][0]{$0.2$}
\psfrag{V}[m][][1][0]{$0.3$} \psfrag{W}[m][][1][0]{$\sigma_{r}$}
\scalebox{1}{\includegraphics{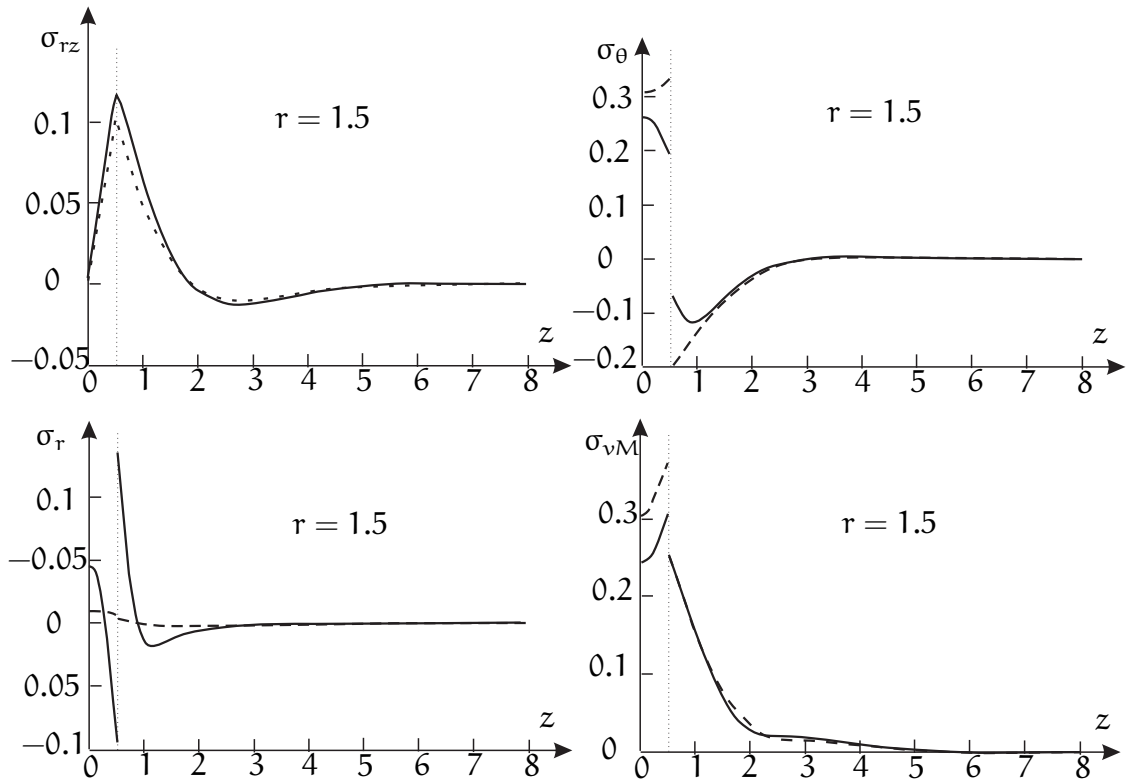}} \caption{Comparison of
the approximate analytical solution (dashed line) with numerical
solution, obtained by the Ritz method (solid line) \label{fig4}}
\end{figure}

We use $r=1,5$ to illustrate on fig. \ref{fig4} the distribution
of stress components along the $z$ axis. Along $r=1,5$ the
$\sigma_{z}$ stress is negligibly small for both solutions, therefore
we do not plot this component.

Hypothesis \eqref{hypothesis} imposes essential restrictions on
the class of solutions, nevertheless the value of shear stress
$\sigma_{rz}$ is captured very well. The error for the hoop stress
$\sigma_{\theta}$ and the radial $\sigma_r$ stress is significant in
the vicinity of the interface $\{z=h\}$.

The approximate analytical solution gives only a smooth part of
the stress field. This solution can be more useful in case of
smooth changing material properties.

\subsection{Solution of the nonlinear problem}
We investigate the error of the perturbation method
$o(s)=\boldsymbol\sigma(s)-(\boldsymbol\sigma^0+ s \boldsymbol\sigma^1)$.

\begin{figure}
\psfrag{P}[m][][1][0]{pressure $p$}
\scalebox{1}{\includegraphics{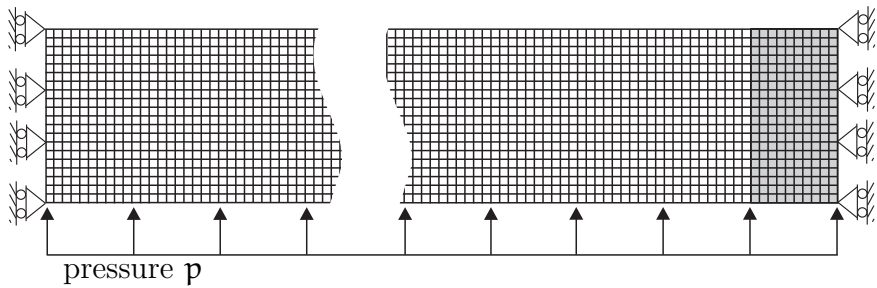}} \caption{FE mesh and boundary conditions \label{fig5}}
\end{figure}

Numerical solutions $\boldsymbol\sigma(s)$ of the system
\eqref{equilInvar} --- \eqref{property} are obtained with the help
of ANSYS finite element code for a set of $s$. The geometry of the half of the
pipe was represented by 3200 axisymmetric PLANE183 finite
elements (fig. \ref{fig5}). We used a uniform mesh with 160 elements along the axis direction z
and 20 elements along the radial direction r.
This type of element can model creep behavior but the
Norton's constitutive law is available only as a secondary creep
equation. That is, the total strain is a sum of elastic strain and
creep strain
\begin{equation*}
\boldsymbol \varepsilon = \boldsymbol \varepsilon^{el} + \boldsymbol \varepsilon^{cr}.
\end{equation*}
If the applied load remains the same with time $t$, then we have
\begin{equation*}
\boldsymbol \varepsilon \approx \boldsymbol \varepsilon^{cr} \ \text{as}  \
t\rightarrow \infty,
\end{equation*}
and we obtain the steady-creep solution as $t\rightarrow \infty$.
In calculations we used $E=100$ for Young modulus and $\nu=0.3$
for Poisson ratio to simulate the elastic material. We consider
the solution to be close enough to the asymptotic solution at the
moment $T=100$ of time $t$.

The solution $\boldsymbol\sigma^1$ of the auxiliary problem is given by
the Ritz method (see section 5.2). Here we used series
\eqref{ser1} --- \eqref{ser2} with $N_r=N_z=25$.

\begin{figure}
\psfrag{A}[m][][1][0]{$0$} \psfrag{B}[m][][1][0]{$1$}
\psfrag{C}[m][][1][0]{$2$} \psfrag{D}[m][][1][0]{$3$}
\psfrag{E}[m][][1][0]{$4$} \psfrag{F}[m][][1][0]{$5$}
\psfrag{G}[m][][1][0]{$6$} \psfrag{H}[m][][1][0]{$7$}
\psfrag{K}[m][][1][0]{$8$} \psfrag{L}[m][][1][0]{$z$}
\psfrag{M}[m][][1][0]{$\sigma_{rz}$}
\psfrag{N}[m][][1][0]{$\sigma_{\theta}$}
\psfrag{O}[m][][1][0]{$s=0.1$} \psfrag{P}[m][][1][0]{$0.005$}
\psfrag{R}[m][][1][0]{$0.01$} \psfrag{S}[m][][1][0]{$-0.005$}
\psfrag{T}[m][][1][0]{$s=0.5$} \psfrag{U}[m][][1][0]{$-0.02$}
\psfrag{X}[m][][1][0]{$0.02$} \psfrag{Y}[m][][1][0]{$0.04$}
\psfrag{Z}[m][][1][0]{$0.06$} \psfrag{V}[m][][1][0]{$s=0.9$}
\psfrag{U2}[m][][1][0]{$-0.05$} \psfrag{X2}[m][][1][0]{$0.05$}
\psfrag{Y2}[m][][1][0]{$0.1$} \psfrag{Z2}[m][][1][0]{$0.15$}
\psfrag{Z3}[m][][1][0]{$0.2$} \psfrag{Z4}[m][][1][0]{$0.25$}
\psfrag{A1}[m][][1][0]{$1.01$} \psfrag{P1}[m][][1][0]{$1.02$}
\psfrag{R1}[m][][1][0]{$1.03$} \psfrag{R2}[m][][1][0]{$1.04$}
\psfrag{S1}[m][][1][0]{$1$} \psfrag{U3}[m][][1][0]{$0.9$}
\psfrag{A3}[m][][1][0]{$0.95$} \psfrag{Y3}[m][][1][0]{$1.05$}
\psfrag{Z3}[m][][1][0]{$1.1$} \psfrag{Z4}[m][][1][0]{$1.15$}
\psfrag{Z5}[m][][1][0]{$1.2$} \psfrag{W1}[m][][1][0]{$0.6$}
\psfrag{W2}[m][][1][0]{$0.8$} \psfrag{W3}[m][][1][0]{$1.0$}
\psfrag{W4}[m][][1][0]{$1.2$} \psfrag{W5}[m][][1][0]{$1.4$}
\psfrag{W6}[m][][1][0]{$1.6$} \psfrag{O5}[m][][1][0]{$r=1.5$}
\psfrag{X4}[m][][1][0]{$0.25$} \psfrag{X3}[m][][1][0]{$0.2$}
\scalebox{0.9}{\includegraphics{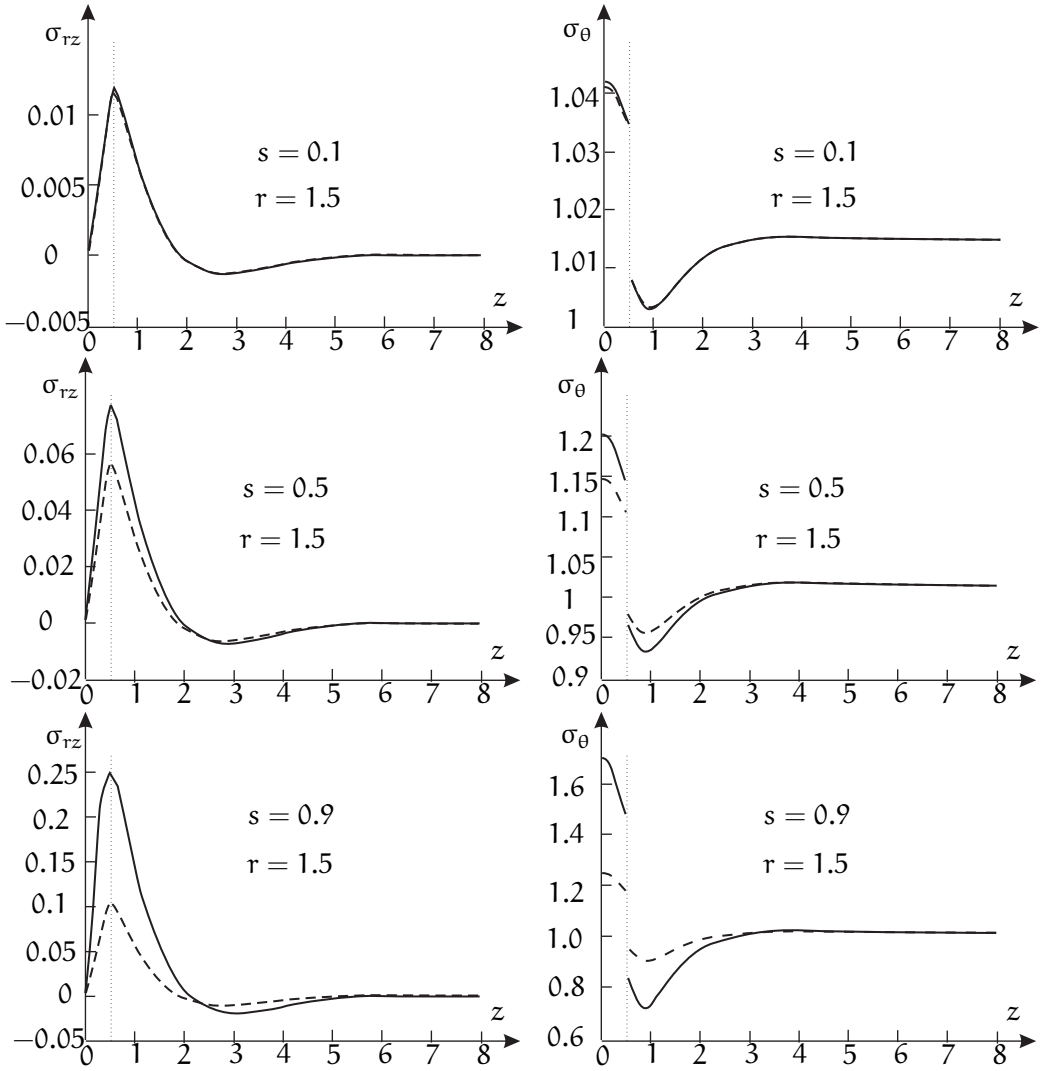}} \caption{Inconsistency
between $\boldsymbol\sigma^0+ s \boldsymbol\sigma^1$ (dashed line) and
$\boldsymbol\sigma(s)$ (solid line) \label{fig6}}
\end{figure}

The leading term in the asymptotic series $\label{hypot}$ is a
good approximation even if the "small" parameter $s$ equals 0.5
(fig. \ref{fig6}). With subsequent increase of $s$ the error
grows dramatically.

\section{Discussion}

The application of the perturbation method to the steady-state
creep problem was investigated. High performance of this method in
prediction of creep response was validated. The perturbation method allow one to
reduce an initial nonlinear problem to the sequence of simpler ones.

This technique is especially attractive if
the unperturbed solution is given in closed form as it was in this paper.
Another example is the creep response of the thick-walled homogeneous pipe
under plane stress conditions \citep{Malinin1981}.
Such solution could be used for perturbation analysis
of creep in open-ended pipes.

The error of the perturbation method becomes substantial when the creep
properties differ from one another by one order of magnitude.
Nevertheless we note that asymptotic expansion \eqref{hypotes}
gives a good simplified model of structure response. This model
treats changes in parameter distribution as jumps of displacements
in linear elastic material.

In that way, the solution for every complicated case of parameter
distribution is represented as a combination of simple solutions.

\section*{Acknowledgments}
A.V. Shutov is grateful for the support provided by the
German Academic Exchange Service.

\appendix
\section*{Appendix}
Assume $\varphi_1(r)=\psi_1(r)=\sin(\frac{r-r_i}{r_o-r_i} \pi)$,  $r_i=1$, $r_o=2$.
After integration in \eqref{eq1}, \eqref{eq2}, \eqref{eq22}
we have the following values of constants
\begin{equation*}
a_1\approx 41.134, \quad a_2\approx 3.770, \quad a_3\approx -0.237,
\end{equation*}
\begin{equation*}
b_1\approx 3.770, \quad b_2\approx -0.237, \quad b_3 \approx 3.948, \quad
b_4 \approx 0.780, \quad b_5 \approx 1.778,
\end{equation*}
\begin{equation*}
k_1 \approx 3.602, \quad k_2 \approx -0.736,\quad  k_3 \approx 1.777,
\end{equation*}
\begin{equation*}
e_2 = 0,\quad  e_3 \approx -1,777, \quad  g_2 \approx 1,777.
\end{equation*}
The ODE \eqref{eq4} has four linearly independent solutions
\begin{equation*}
\psi^{I}_2= e^{Re(\lambda) z} \sin(Im(\lambda) z), \quad \psi^{II}_2= e^{Re(\lambda) z} \cos(Im(\lambda) z),
\end{equation*}
\begin{equation*}
\psi^{III}_2= e^{-Re(\lambda) z} \sin(Im(\lambda) z), \quad \psi^{IV}_2= e^{-Re(\lambda) z} \cos(Im(\lambda) z),
\end{equation*}
where $\lambda$ is one of solutions of characteristic equation
\begin{equation*}
k_3 \lambda^4 + k_2 \lambda^2 + k_1 = 0.
\end{equation*}
To be definite, we use
\begin{equation*}
\lambda = 0.903 + 0.780 \ i.
\end{equation*}

\address{Alexey Shutov,
Lavrentyev Institute of Hydrodynamics of SB RAS, 630090, Novosibirsk, Russia;
Konstantin Naumenko, Holm Altenbach, Lehrstuhl f\"ur Technische Mechanik
Fachbereich Ingenieurwissenschaften Martin-Luther-Universit\"at
Halle-Wittenberg D-06099 Halle}

\end{document}